\documentclass{amsart}
\usepackage{amssymb}
\usepackage{amsfonts}
\usepackage{latexsym}
\usepackage[all]{xy}
\usepackage{verbatim}

\newtheorem{theorem}{Theorem}[section]
\newtheorem{lemma}[theorem]{Lemma}
\newtheorem{proposition}[theorem]{Proposition}
\newtheorem{corollary}[theorem]{Corollary} 
\theoremstyle{definition}  
\newtheorem{definition}[theorem]{Definition}

\newtheorem{remark}[theorem]{Remark}

\newcommand{\Tr}{\text{Tr}}
\newcommand{\id}{\text{id}}

\newcommand{\End}{\text{End}}

\newcommand{\Hom}{\text{Hom}} 
\newcommand{\uHom}{\underline{\text{Hom}}}

\newcommand{\Rep}{\text{Rep}}

\newcommand{\op}{\text{op}}
 
\newcommand{\dH}{\dim \Hom} 

\newcommand{\ev}{\text{ev}}
\newcommand{\coev}{\text{coev}}

\newcommand{\C}{\mathcal{C}}
\newcommand{\D}{\mathcal{D}}
\newcommand{\Z}{\mathcal{Z}}
\renewcommand{\H}{\mathcal{H}}
\newcommand{\M}{\mathcal{M}}

\newcommand{\be}{\mathbf{1}}

\newcommand{\actl}{\rightharpoonup}
\newcommand{\actr}{\leftharpoonup}

\renewcommand{\1}{_{(1)}} 
\renewcommand{\2}{_{(2)}}

\newcommand{\mC}{{\mathcal C}}
\newcommand{\cC}{{\mathcal C}}

\newcommand{\mM}{{\mathcal M}}
\begin{document}
\title[A categorical analogue of Radford's $\mathbf{S^4}$-formula]
{An analogue of  Radford's $\mathbf{S^4}$-formula for finite tensor categories}

\author{Pavel Etingof}
\address{Department of Mathematics, Massachusetts Institute of Technology,
Cambridge, MA 02139, USA}
\email{etingof@math.mit.edu}

\author{Dmitri Nikshych}
\address{Department of Mathematics and Statistics,
University of New Hampshire,  Durham, NH 03824, USA}
\email{nikshych@math.unh.edu}

\author{Viktor Ostrik}
\address{Institute for Advanced Study,
Einstein Dr, Princeton, NJ 08540, USA}
\email{ostrik@math.mit.edu}

\date{April 27, 2004}
\begin{abstract}
We develop the theory of Hopf bimodules for a finite rigid tensor category $\C$.
Then we use this theory to define a distinguished invertible
object $D$ of $\C$ and an isomorphism of tensor functors $\delta:
V^{**}\to D\otimes ^{**}V\otimes D^{-1}$. This provides a
categorical generalization of Radford's $S^4$ formula 
for finite dimensional Hopf algebras \cite{R1},
which was proved in \cite{N} for weak Hopf
algebras, in \cite{HN} for quasi-Hopf algebras, 
and conjectured in general in \cite{EO}. 
When $\C$ is braided, we establish a connection between 
$\delta$ and the Drinfeld isomorphism of $\C$, extending the result of \cite{R2}.
We also show that a factorizable braided tensor category is
unimodular (i.e. $D={\mathbf 1}$).
Finally, we apply our theory to prove that the pivotalization of a fusion
category is spherical, and give a purely  algebraic characterization of exact module 
categories defined in \cite{EO}.
\end{abstract} 
\maketitle  

\begin{section}{Introduction}

One of the most important general results about finite
dimensional Hopf algebras is Radford's formula for the forth 
power of the antipode. This formula reads: for any finite dimensional Hopf algebra $H$
over a field $k$ with antipode $S$, 
\begin{equation}
S^4(h) = a^{-1}(\alpha\actl h \actr \alpha^{-1})a, 
\end{equation}
where $a$ and $\alpha$ are the distinguished grouplike elements
of $H$ and $H^*$, respectively. 

Recently, several generalizations of this formula were
discovered. Namely, in \cite{N} the formula was extended to 
weak Hopf algebras, and in \cite{HN} to quasi-Hopf algebras. 
Finally, in \cite{EO} the authors conjectured a
generalization of Radford's formula to any finite tensor category
(Conjecture 2.15). One of the main achievements of the present
paper is a proof of this conjecture. 
\footnote{We note that this result is, essentially, equivalent
to the statement that Radford's formula holds for weak quasi-Hopf
algebras. However, weak quasi-Hopf algebras are very cumbersome
objects, and we avoid working with them by using a much simpler purely
categorical language.}

More specifically, we show in Theorem~\ref{main} (which is our main theorem)
that for a finite tensor category $\C$ there is an invertible 
object $D$, called the {\em distinguished invertible object} of
$\C$, and an isomorphism 
of tensor functors $\delta :\, ?^{**} \cong D\otimes {}^{**}? \otimes
D^{-1}$. We also show that our definition of the distinguished
object is the same as in \cite[Definition~2.12]{EO}, which means
that our main theorem implies Conjecture 2.15 from \cite{EO}. 

Further, we apply the main theorem
to braided tensor categories to show that 
a factorizable braided category is unimodular (in particular, the center
$\Z(\C)$ of a finite tensor category $\C$ is unimodular). We also
relate the isomorphism $\delta$
with the Drinfeld isomorphism $ u :\, ?\cong ?^{**}$. 

In the case when $\C$ is semisimple,
we describe the isomorphism $\delta$ in terms
of M\"uger's squared norms of simple objects,  
and show that the pivotalization
of $\C$ is spherical (this was previously known only in 
characteristic $0$ \cite{ENO}).

The paper is organized as follows.  In Section~\ref{section:Hopf bimodules} we define 
a category $\H$ of Hopf bimodules over $\C$ as the category of right modules
over the algebra $A= \uHom({\mathbf 1},\,{\mathbf 1})$ in $\C \boxtimes \C^{\op}$. 
We define a tensor category structure on $\H$ and prove in
Proposition~\ref{equivalence} that there is a tensor equivalence between $\C$ and $\H$ 
(which establishes a categorical version of the Fundamental Theorem for Hopf bimodules).
We use this result as a major tool in our arguments throughout the paper.

In Section~\ref{construction of the isomorphism} we define the 
{\em distinguished invertible object} of $\C$ as the unique up to an isomorphism
object $D$ such that $A^* \cong (D\boxtimes {\mathbf 1})\otimes A$, following
the idea used in \cite{HN} for quasi-Hopf algebras. We use the category of Hopf
bimodules to construct in Theorem~\ref{main} a natural tensor
isomorphism $\delta :\, ?^{**} \cong D\otimes \,{}^{**}? \otimes D^{-1}$.

In Sections~\ref{section:unimodularity} and \ref{Drinfeld iso} we work with a braided
finite tensor category $\C$. We introduce a categorical notion of {\em factorizability}
of $\C$ and show that the center of a finite tensor category is factorizable. We prove
that a factorizable tensor category is automatically unimodular (i.e., $D\cong {\mathbf 1}$),
extending the results known for (weak, quasi) Hopf algebras. We also show that
in the case of unimodular $\C$ one has $\delta = u^{-1}\circ u^*$ where $u :\, 
?\cong \,?^{**}$
is the Drinfeld isomorphism.

In Section~\ref{comp EO} we check that our definition of $D$ agrees with 
the one given in \cite{EO}, that is, the projective cover of the unit object
${\mathbf 1}$ coincides with the injective hull of~$D$.

In Section~\ref{fusion categories} we specialize to the case when $\C$ is semisimple
and give a convenient numerical characterization of the isomorphism $\delta$ in
Theorem~\ref{vitia monstr}. As a consequence, we obtain that the pivotalization 
of a semisimple category is spherical.

The appendix contains an algebraic characterization of exact module categories
studied in \cite{EO}. We show in Theorem~\ref{M exact iff} that if $\M$ is a module category
over $\C$, $B$ is a $k$-algebra such that $\M\cong B-\mbox{mod}$, 
and $\bar{F} :\C \to\mbox{Vect}_k$  is the fiber functor  constructed
from the action of $\C$ on $\M$, then $\M$ is exact if and only if 
the algebra $H = \End(\bar{F})$ is a projective $B\otimes B^\circ$-module.

{\bf Remark.} As was anticipated by G. Kuperberg, 
the statement and proof of the categorical Radford's formula do not
make an essential use of the additive structure of the category
$\C$ (and thus can be generalized to the non-additive case). 
In this respect, they are similar to Theorem 3.10 in
\cite{Ku}, which gives Radford's formula for a Hopf object in a
rigid monoidal (possibly non-additive) category. 

\textbf{Acknowledgements.} 
We are grateful to G. Kuperberg for a useful discussion.
The research of P.E. was supported by the NSF
grant DMS-9988796.
The research of D.N. was supported by the NSF
grant DMS-0200202.
The research of V.O. was 
supported by the NSF grants DMS-0098830 and DMS-0111298.
\end{section}

\begin{section}{An analogue of the category of Hopf bimodules}
\label{section:Hopf bimodules}

We work over an algebraically closed field $k$.  All categories considered 
in this paper are abelian over $k$ with all objects being of finite length 
and all morphism spaces being finite dimensional.  

Such a category $\C$ is said to be {\em finite} 
if it has finitely many isomorphism classes of simple objects and every object
has a projective cover. That is, $\C$ is equivalent to the category 
of finite dimensional representations of some finite dimensional $k$-algebra.

By a {\em tensor category} over $k$ we understand an abelian rigid monoidal 
category. We refer the reader to \cite{BK, K, EO}  for the general theory of 
tensor categories and module categories over them.

Recall the notion of Deligne's
tensor product $\C \boxtimes \D$ of abelian categories $\C,\D$, see \cite{D}. 
By definition, this is the universal object for the functor assigning to every
abelian category $\mathcal{A}$ the category 
of additive right (equivalently, left) exact bifunctors from $\C\times \D$
to $\mathcal{A}$. If $\C,\D$ are tensor categories 
then the category $\C \boxtimes \D$ has a natural structure of 
a   tensor category with the tensor product  
\begin{equation}
\otimes \times \otimes : (\C \boxtimes \D)\times  (\C \boxtimes \D)
\to \C \boxtimes \D,
\end{equation}
(which we will still denote $\otimes$),  and with the unit object
${\mathbf 1}\boxtimes {\mathbf 1}$.

Deligne's tensor product can also be applied to functors. 
Namely, if $F: \C\to \C'$ and $G: \D\to \D'$ are additive right (left) exact
functors between abelian categories then one can 
define the functor $F\boxtimes G: \C\boxtimes D\to 
\C'\boxtimes \D'$. 

Let $\C$ be a   finite tensor category over $k$. Let 
$\C^{\op}$ be the opposite tensor category, that is,
$\C^{\op}=\C$ as a category but the tensor product $\otimes^{\op}$ in $\C^{\op}$
is different: $X\otimes^{\op} Y=Y\otimes X$. 

The category $\C$ has a natural structure of an exact module category \cite{EO} over
$(\C \boxtimes \C^{\op},\, \otimes)$ coming from
\begin{equation}
\otimes\circ(\otimes \times\id_{\C^{op}}):\C\times\C\times\C^{op}\to \C.
\end{equation}

We will denote this action by $(X,V)\mapsto X\circ V$ ($X\in
\C\boxtimes \C^{op}$, $V\in \C$). 

The internal Hom, $\uHom(V_1,\, V_2)$,
of two objects $V_1,\, V_2$ in $\C$ is an object of $\C \boxtimes \C^{\op}$
representing the functor $\Hom_\C(?\circ V_1,\, V_2)$  (the latter contravariant
functor is left exact,
so it is representable).  This means that there is a natural isomorphism
\begin{equation}
\Hom_\C(X\circ V_1,\, V_2) \cong \Hom_{\C \boxtimes \C^{\op}}(X,\, \uHom(V_1,\, V_2)).
\end{equation}

We refer to \cite{EO} for the properties of $\uHom$. The object 
$A:=\uHom({\mathbf 1},{\mathbf 1})$ has a
natural structure of an algebra in the category $\C \boxtimes \C^{\op}$. 

\begin{definition} 
\label{Hopf bimodules}
The category of right $A$-modules in  
$(\C \boxtimes \C^{\op},\, \otimes)$
will be called the category of 
{\em Hopf bimodules} over $\C$.
\end{definition}

Let $\H$ denote the category of Hopf bimodules.

Observe that $(X,Y)\mapsto\Hom_\C({}^*Y ,X)$ 
is an additive bifunctor from $\C^{\op} \times \C$
to $\mbox{Vect}_k$. Hence it defines an additive functor 
$H_{\C}: \C^{\op} \boxtimes \C\to \mbox{Vect}_k$.
Therefore, one can define another tensor product $\odot$ on 
$\C \boxtimes \C^{\op}$ by
\begin{equation}
(\id_{\C}\boxtimes H_\C)\boxtimes \id_{\C^{\op}} :
(\C \boxtimes \C^{\op})\boxtimes (\C \boxtimes \C^{\op}) \cong
(\C \boxtimes (\C^{\op}\boxtimes \C) )\boxtimes \C^{\op} \to \C \boxtimes \C^{\op},
\end{equation}
where we implicitly used a natural action of $\mbox{Vect}_k$ on $\C$.
One checks directly that $A$ is the unit object for $\odot$.

\begin{remark} \label{functors}
(i) Another way to define $\odot$ is the following: one identifies $\C^{op}$
with the dual category $\C^\vee$ via the functor $X\mapsto {}^*X$. Then the
category $\C \boxtimes \C^{op}$ is identified with the category of left exact
functors from $\C$ to itself (for example $X\boxtimes Y\in 
\C \boxtimes C^{op}$ corresponds to the functor $Z\mapsto 
\Hom ({}^*Z,Y)\otimes X$), and $\odot$ corresponds to the composition
of functors. Under this identification the object $A\in \C \boxtimes \C^{op}$
representing the functor $X\boxtimes Y\mapsto \Hom (X\otimes Y, \be)=
\Hom (X,Y^*)$ corresponds to the identity functor (this is why we use 
$X\mapsto {}^*X$ and not $X\mapsto X^*$ to identify $\C^{op}$ and $\C^\vee$).

(ii) If $\C$ is the representation category of a Hopf algebra 
(or quasi-Hopf algebra, or weak Hopf algebra) $H$
then $\H$ is the  category of usual $H$-Hopf bimodules \cite{LS,
HN, BNS} and $\odot$ is dual to the usual bimodule tensor product
(thus $M\odot N=(N^*\otimes_HM^*)^*$ where the star denotes the dual
vector space).
\end{remark}

The Fundamental Theorem for Hopf modules over a Hopf algebra $H$ \cite{LS} states
that the category of $H$-Hopf modules is equivalent to the category of vector
spaces. In \cite{HN} it was explained that for quasi-Hopf algebras the notion of
a Hopf module should be replaced by that of a Hopf {\em bimodule} and the Fundamental
Theorem for Hopf bimodules over a quasi-Hopf algebra \cite[Proposition 3.11]{HN} 
algebra was proved. Proposition~\ref{equivalence} below provides a categorical
version of the Fundamental Theorem for Hopf bimodules and
generalizes 
the results of 
\cite{LS, HN}.

\begin{proposition}
\label{equivalence} 
(a) The category $\H$ is a tensor subcategory of $(\C \boxtimes \C^{\op},\, \odot)$ 
and the functors $(?\boxtimes {\mathbf 1})\otimes A$ and $({\mathbf 1}\boxtimes ?)\otimes A$
from $\C$ to $\H$  are equivalences of  tensor categories. 
\newline\qquad (b) There is a natural isomorphism of tensor functors  
\begin{equation}
\label{RV}
\rho : (?\boxtimes {\mathbf 1})\otimes A \cong ({\mathbf 1}\boxtimes ?)\otimes A.
\end{equation}
\end{proposition}
\begin{proof} 
(a) For all objects $V$ in $\C$ we have
\begin{eqnarray*}
\uHom({\mathbf 1},V)
&=& \uHom({\mathbf 1},(V\boxtimes {\mathbf 1})\circ {\mathbf 1})
\cong (V\boxtimes {\mathbf 1})\otimes \uHom({\mathbf 1},{\mathbf 1})
= (V\boxtimes {\mathbf 1})\otimes A, \\
\uHom({\mathbf 1},V)
&=& \uHom({\mathbf 1},({\mathbf 1}\boxtimes V) \circ {\mathbf 1})
\cong ({\mathbf 1}\boxtimes V)\circ \uHom({\mathbf 1},{\mathbf 1})
= ({\mathbf 1}\boxtimes V)\otimes A.
\end{eqnarray*}
It follows from \cite{EO}, Theorem 3.17 that $(?\boxtimes {\mathbf 1})\otimes A$
is an equivalence between $\C$ and $\H$. 
To see that it is tensor, observe that 
under the identification in Remark \ref{functors} (i) the functor $(?\boxtimes
\be)\otimes A$ (and similarly $(\be \boxtimes ?)\otimes A$) sends $V\in \C$ to
the functor $V\otimes ?$ from $\C$ to itself. Thus the associativity constraint
in the category $\C$ gives rise to a tensor structure on these functors.


(b) The tensoriality of the natural isomorphism \eqref{RV} is obvious from
description in (a). It is also equivalent to 
commutativity of the following diagram
\begin{equation}
\xymatrix{ ((V\boxtimes {\mathbf 1})\otimes A)\odot ((W\boxtimes {\mathbf 1})\otimes A) 
\ar[r]
\ar[d]^{\rho_V\odot\rho_W}&
((V\otimes W)\boxtimes {\mathbf 1})\otimes A \ar[d]^{\rho_{V\otimes W}} \\
(({\mathbf 1}\boxtimes V)\otimes A)\odot (({\mathbf 1}\boxtimes W)\otimes A) \ar[r]& 
 ({\mathbf 1}\boxtimes  (V\otimes W))\otimes A. }
\end{equation}
\end{proof}

\begin{remark} Another way to state Proposition \ref{equivalence} (b) is to
say that the following diagram commutes:
\begin{equation} \label{tens}
\xymatrix{ (V\boxtimes {\mathbf 1})\otimes (W\boxtimes {\mathbf 1})\otimes A 
\ar[r]
\ar[d]^{\id \otimes \rho_W}&
((V\otimes W)\boxtimes {\mathbf 1})\otimes A \ar[d]^{\rho_{V\otimes W}} \\
(V\boxtimes \be)\otimes (\be \boxtimes W)\otimes A \ar[d]&
(\be \boxtimes (V\otimes W))\otimes A  \\
(\be \boxtimes W)\otimes (V\boxtimes \be)\otimes A \ar[r]^{\id \otimes 
\rho_V} &
(\be \boxtimes W)\otimes (\be \boxtimes V)\otimes A. \ar[u] } 
\end{equation}
\end{remark}

%

\end{section}
\begin{section}
{Construction of an isomorphism between duality functors}
\label{construction of the isomorphism}

Let $A=\uHom({\mathbf 1},\,{\mathbf 1})$ be the algebra in $\C \boxtimes \C^{\op}$ 
defined in the previous 
Section and let $M$ be a left $A$-module. Recall \cite{O1} that $M^*$ has a natural structure
of a right $A-$module with the action given by
\begin{equation}
\label{dual module}
M^* \otimes A  \stackrel{m_A^*\otimes \id_A}{\longrightarrow} M^* \otimes A^*\otimes A
\stackrel{id_{M^*}\otimes \coev_A}{\longrightarrow} M^*,
\end{equation}
where $m_A : A\otimes M \to M$ be the left action of $A$ on $M$ and
$\coev_A$ is the coevaluation morphism of $A$.
In particular $A^*$ has a canonical
\footnote{In this paper, ``canonical''  
means that there is a distinguished 
choice which should be obvious to the reader.} 
structure of a Hopf bimodule. Thus according to Proposition~\ref{equivalence} (a)
there exists a unique up to an isomorphism object $D\in \C$ such that 
\begin{equation}
\label{Definition of D}
(D\boxtimes {\mathbf 1})\otimes A \cong A^*
\end{equation}
as Hopf bimodules. Moreover isomorphism \eqref{Definition of D} is unique up
to scaling. It follows immediately from the definition that the
Frobenius-Perron dimension (see \cite{EO}) of $D$ equals to 1 and thus $D$
is invertible. 

\begin{definition} 
\label{dist inv o}
The object $D$ defined by \eqref{Definition of D}
is called a {\em distinguished invertible object} of $\C$.
\end{definition}

\begin{remark}
Our definition of the distinguished invertible object of $\C$ categorically
extends definitions of distinguished group-like element, or modulus, of a
Hopf algebra \cite{R1}, quasi-Hopf algebra \cite{HN}, or weak Hopf algebra 
\cite{N}. In Section~\ref{comp EO} we show that Definition~\ref{dist inv o} 
is equivalent to \cite[Definition 2.12]{EO}. 
\end{remark}

The classical formula of D.~Radford \cite{R1} expresses the fourth power of the antipode 
of a finite-dimensional Hopf algebra $H$ in terms of distinguished group-like elements
of $H$ and $H^*$.  The categorical version of this formula below is a main result of this
paper.

\begin{theorem}
\label{main}
Let $\C$ be a   finite tensor category. There is a natural isomorphism
of tensor functors, 
\begin{equation}
\label{dV}
\delta\,  :\, ?^{**}\to D\otimes {}^{**}?\otimes D^{-1}.
\end{equation}
\end{theorem}
\begin{proof}
Isomorphism \eqref{Definition of D} produces a canonical isomorphism of 
algebras 
\begin{equation}
\label{twostar}
A^{**}=\uHom(A^*,A^*)\cong \uHom((D\boxtimes \be)\otimes A, 
(D\boxtimes \be)\otimes A)=(D\boxtimes \be)\otimes A\otimes (D\boxtimes \be)^*.
\end{equation}
We will identify these algebras using this isomorphism.

Recall that we have a tensor isomorphism \eqref{RV} $\rho_V: (V\boxtimes \be)
\otimes A\cong (\be \boxtimes V)\otimes A$. Its double dual $\rho_V^{**}:
(V^{**}\boxtimes \be)\otimes A^{**}\cong (\be \boxtimes {}^{**}V)\otimes 
A^{**}$ is also tensor (i.e. the diagram analogous to \eqref{tens} commutes).

Thus we have a tensor isomorphism of right $A-$modules
$$\tilde \rho_V: (V^{**}\boxtimes \be)\otimes (D\boxtimes \be)\otimes A \cong
(\be \boxtimes {}^{**}V)\otimes (D\boxtimes \be)\otimes A$$
defined by $\tilde \rho_V\otimes \id_{(D\boxtimes \be)^*}=\rho_V^{**}$.

Now define $\tilde \delta_V$ as the following composition:
$$(V^{**}\boxtimes \be)\otimes A=((V^{**}\otimes D^*)\boxtimes \be)\otimes 
(D\boxtimes \be)\otimes A\stackrel{\tilde \rho_{V\otimes D^*}}{\longrightarrow}
(\be \boxtimes ({}^{**}V\otimes D^*))\otimes (D\boxtimes \be)\otimes A=$$
$$=(D\boxtimes \be)\otimes (\be \boxtimes ({}^{**}V\otimes D^*))\otimes A
\stackrel{\id \otimes \rho_{{}^{**}V\otimes D^*}}{\longrightarrow}
(D\boxtimes \be)\otimes (({}^{**}V\otimes D^*)\boxtimes \be)\otimes A=$$
$$=((D\otimes {}^{**}V\otimes D^*)\boxtimes \be)\otimes A.$$ 
Obviously, the isomorphism $\tilde \delta_V$ is tensor (again, the diagram
analogous to \eqref{tens} commutes). 

Finally, define the isomorphism $\delta_V: V^{**}\to D\otimes {}^{**}V
\otimes D^*$ by the condition $\delta_V\otimes \id_A=\tilde \delta_V$.
Since $\tilde \delta_V$ is a morphism of right $A-$modules Proposition
\ref{equivalence} (a) implies that $\delta_V$ is well defined. The fact that
$\tilde \delta_V$ is tensor translates into the fact that $\delta_V$ is an
isomorphism of tensor functors. The Theorem is proved.



\end{proof}

\begin{corollary}
There is a positive integer $N$ such that the $N$th powers of 
tensor functors  $?^{**}$ and ${}^{**}?$ are naturally isomorphic.
\end{corollary}
\begin{proof}
Since $\C$ has finitely many non-isomorphic invertible objects, there exists
$N$ such that $D^{\otimes N} \cong {\mathbf 1}$.
\end{proof}

\end{section}


\begin{section}
{Unimodularity of factorizable categories}
\label{section:unimodularity}

\begin{definition}
\label{unimodularity}
We will say that $\C$ is {\em unimodular} if its distinguished invertible
object is isomorphic to ${\mathbf 1}$.
\end{definition}

Equivalently, $\C$ is unimodular if $A$ is a self-dual object of $\C\boxtimes \C^{\op}$.

Let $(\C,\, \sigma)$ be a braided finite tensor category, where $\sigma$ is a 
natural isomorphism of bifunctors $\sigma :\otimes\cong \otimes^{\op}$ 
satisfying hexagon axioms \cite{BK, K}.  Let $\Z(\C)$ be the center of $\C$.
Recall (see e.g. \cite{K}) that the objects of $\Z(\C)$ are pairs $(X,e_X(?))$
where $e_X(?)$ is a functorial isomorphism $e_X(Y): X\otimes Y\simeq Y\otimes X$ defined
for all $Y\in \C$ and satisfying certain axioms. Now we define a tensor
functor $G: \C \boxtimes \C^{op} \to \Z(\C)$ in the following way:
\begin{equation}
\label{G}
G(X\boxtimes Y)=(X\otimes Y,e_{X\otimes Y})
\end{equation}
 where
$$
e_{X\otimes Y}(Z): X\otimes Y\otimes Z\stackrel{\id_X\otimes 
\sigma_{Z,Y}^{-1}}{\longrightarrow} X\otimes Z\otimes Y \stackrel{\sigma_{X,Z}
\otimes \id_Y}{\longrightarrow} Z\otimes X\otimes Y.$$
The functor $G$ has a natural structure of a braided tensor functor.


\begin{definition}
\label{factorizability}
We will say that a finite braided tensor category $\C$ is {\em factorizable} if
$G$ is an equivalence of tensor categories.
\end{definition}

\begin{remark} 
Factorizable Hopf algebras were introduced and studied by N.~Reshetikhin and
M.~Semenov-Tian-Shansky in \cite{RS}. This notion was extended to weak Hopf
algebras in \cite{NTV} and to quasi-Hopf algebras in \cite{BT}. One can directly
check that our Definition~\ref{factorizability} extends the previous definitions.
E.g., using a computation analogous to the one
given by H.-J.~Schneider for Hopf algebras in \cite[Theorem 4.3]{S}
one shows that a weak Hopf algebra $H$ is factorizable if and only if
its representation category $\Rep(H)$ is factorizable in the sense of 
Definition~\ref{factorizability}, so the two definitions agree in this case.

The notion of a factorizable braided tensor category 
also extends that of a modular category to the case when $\C$ is not necessarily
ribbon or semisimple. Indeed, every  semisimple finite tensor category 
is equivalent to the representation category of some semisimple weak Hopf algebra 
\cite{O1}, and it was shown in \cite{NTV} that a semisimple ribbon weak Hopf algebra
is modular if and only if it is factorizable.
\end{remark}

An example of a factorizable category is given by the center of a tensor
category.

\begin{proposition}
\label{Z is factorizable}
Let $\C$ be a  finite tensor category. Then its center $\Z(\C)$ is factorizable.
\end{proposition}
\begin{proof}
Let $\M$ be an exact
module category over a  finite tensor category
$\D$ and let $\D^*_\M$ be the dual tensor category,
whose objects are $\D$-module endofunctors of $\M$, see \cite{EO, O1} for definitions.
By \cite[Corollary 3.35]{EO} there is a canonical tensor equivalence 
of tensor categories $ Q: \Z(\D)\cong \Z(\D^*_\M)$
that assigns to every object in $\Z(\D)$ its module action on $\M$. 
We will apply this result to the case when $\D = \C\boxtimes \C^{\op}$ and $\M=\C$.

Recall that  there is a tensor equivalence 
$\Z(\C)\cong (\C\boxtimes \C^{\op})^*_\C$ \cite{EO, O1}. 
It is straightforward to check that the tensor 
functor $G : \Z(\C) \boxtimes \Z(\C)^\op \to \Z(\Z(\C))$ defined as in 
equation \eqref{G} is the composition of the  obvious equivalence 
$\Z(\C)\boxtimes \Z(\C)^{\op} \cong  \Z(\C\boxtimes \C^{\op})$ 
and $Q : \Z(\C\boxtimes \C^{\op}) \cong \Z(\Z(\C))$ defined in the previous
paragraph. Therefore, $G$ is an equivalence, i.e., $\Z(\C)$ is factorizable.
\end{proof}

Next, we establish a categorical generalization of another Radford's
result \cite{R3} stating that a factorizable Hopf algebra is unimodular
(see also \cite{BT} for quasi-Hopf algebras).

\begin{proposition}
If $\C$ is a factorizable tensor category then $\C$ is unimodular.
\end{proposition}
\begin{proof}
Observe that $\C$ is 
a $\Z(\C)$-module category via the forgetful functor $F: \Z(\C)\to \C$,
\begin{equation}  
\label{forgetful}
Z \bullet V := F(Z)\otimes V,
\end{equation}
for all objects $Z\in \Z(\C)$ and $V$ in $\C$. By the factorizability of $\C$
there is a natural isomorphism of functors $?\bullet {\mathbf 1} \cong G^{-1}(?)\circ {\mathbf 1}$.
Let $I: C\to \Z(\C)$ be the induction functor \cite[Lemma 3.38]{EO}, i.e., the
right adjoint functor of $F$. There is a sequence of natural isomorphisms :
\begin{eqnarray*}
\Hom_{\Z(\C)}(Z,\, I(V))
&\cong& \Hom_{\C}(F(Z),\, V)\\
&=& \Hom_{\C}(Z\bullet {\mathbf 1},\, V)\\
&\cong& \Hom_{\C}(G^{-1}(Z)\circ {\mathbf 1},\, V)\\
&\cong& \Hom_{\C\boxtimes \C^{\op}}(G^{-1}(Z),\, \uHom({\mathbf 1},\,V))\\
&\cong&  \Hom_{\Z(\C)}(Z,\,G(\uHom({\mathbf 1},\,V))).
\end{eqnarray*}
Hence, $I(?)$ is naturally isomorphic to $G(\uHom({\mathbf 1},\,?))$. 
For all objects $V$ in $\C$ we have $F({}^*V) = {}^*F(V)$ 
by the definition of duality in $\Z(\C)$, i.e., $F$ commutes
with the left dual functor. Therefore the adjoint functor $I$
commutes with the right dual functor. In particular, $I({\mathbf 1})^*= I({\mathbf 1})$
and $\uHom({\mathbf 1},\,{\mathbf 1}) = \uHom({\mathbf 1},\,{\mathbf 1})^*$ (note that the tensor functor $G$ 
commutes with duality), i.e., $\C$ is unimodular.
\end{proof}

\end{section}

\begin{section}{Relation with the Drinfeld isomorphism in the braided case}
\label{Drinfeld iso}

Let $\C$ be a braided tensor category with braiding $\sigma : \otimes \cong \otimes^{\op}$. 
It is well known that in this case
there is a natural isomorphism $u: \,? \to ?^{**}$, called the {\em Drinfeld isomorphism},
given by 
\begin{equation}
\label{uV}
u_V : V  \stackrel{\coev_{V^*}}{\longrightarrow} V \otimes V^*\otimes V^{**}
\stackrel{\sigma_{V,V^*}}{\longrightarrow} V^* \otimes V \otimes  V^{**} 
\stackrel{\ev_{V}}{\longrightarrow}  V^{**}, 
\end{equation}
where $\coev_V : {\mathbf 1}\to V\otimes V^*$ and $\ev_V :  V^* \otimes V \to {\mathbf 1}$
are the coevaluation and evaluation morphisms attached to an object $V$.

Let $\C =\Rep(H)$ be the category of representations of a finite-dimensional
Hopf algebra $H$. Then a result of Radford \cite{R2} relates the Drinfeld
isomorphism $u_V : V\cong V^{**}$ and the tensor isomorphism $\delta_V : V^{**}
\to D \otimes {}^{**}V\otimes D^{-1}$ 
constructed in Theorem~\ref{main}. An extension of this result
to weak Hopf algebras was obtained in \cite[Lemma 5.12]{ENO}. The proofs 
use Hopf algebra language and techniques. Below we give
a categorical generalization of these results (we restrict ourselves 
to the unimodular case). 

\begin{theorem}
Let $\C$ be a unimodular braided finite tensor category. Let  
$\delta_V : V^{**}\to {}^{**}V$ be the tensor isomorphism 
constructed in Theorem~\ref{main} and let $u_V : V\to V^{**}$ be the Drinfeld
isomorphism. Then there is an equality of natural isomorphisms  
\begin{equation}
\delta_V =  u_{{}^{**}V}^{-1} \circ  u_{{}^*V}^*.
\end{equation}
\end{theorem}
\begin{proof}
Let $\rho_V : (V\boxtimes {\mathbf 1})\otimes A\cong ({\mathbf 1}\boxtimes V)\otimes A$ be as
in Theorem~\ref{main}.
Recall that by definition,   
\begin{equation}
\label{unimodular delight}
(\delta_V\boxtimes \id_{\mathbf 1})\otimes \id_A = \rho_{{}^{**}V}^{-1}\circ {\rho_V}^{**}.
\end{equation}
Let $\Sigma = \sigma\boxtimes \sigma^{-1}$ be the braiding on $\C\boxtimes \C^{\op}$.
Define a natural  isomorphism 
\begin{equation}
\iota : (?\boxtimes {\mathbf 1})\otimes A \to (?^{**}\boxtimes {\mathbf 1})\otimes A
\end{equation}
as the following composition :
\begin{equation}
\label{IV}
\iota_V : (V\boxtimes {\mathbf 1})\otimes A \stackrel{\rho_V}{\longrightarrow}
({\mathbf 1} \boxtimes V)\otimes A  \stackrel{\Sigma_{{\mathbf 1} \boxtimes V, A}}{\longrightarrow}
A \otimes ({\mathbf 1} \boxtimes V)  \stackrel{\rho_{V^*}^*}{\longrightarrow}
A \otimes (V^{**} \boxtimes {\mathbf 1})  \stackrel{\Sigma_{A, V^{**} \boxtimes {\mathbf 1}}}{\longrightarrow}
(V^{**} \boxtimes {\mathbf 1})\otimes A.
\end{equation}
Observe that in terms of $\Hom$ spaces the isomorphism $\iota_V$ is given
as the following sequence of natural isomorphisms: 
\begin{eqnarray*}
\Hom(X_1\boxtimes X_2,\, (V\boxtimes {\mathbf 1})\otimes A )
&\cong& \Hom(V^*\otimes X_1\otimes X_2,\, {\mathbf 1})\\
&\cong& \Hom(X_1\otimes V^* \otimes X_2,\, {\mathbf 1})\\
&\cong& \Hom(X_1\otimes X_2 \otimes V^*,\, {\mathbf 1})\\
&\cong& \Hom(X_1\boxtimes X_2,\, (V^{**}\boxtimes {\mathbf 1})\otimes A )
\end{eqnarray*}
for all objects $X_1, X_2\in \C$, where the two isomorphisms in the middle
come from the braiding in $\C$.  On the other hand the isomorphism
$(u_V\boxtimes \id_{\mathbf 1})\otimes \id_A$ is given by 
\begin{eqnarray*}
\Hom(X_1\boxtimes X_2,\, (V\boxtimes {\mathbf 1})\otimes A )
&\cong& \Hom(V^*\otimes X_1\otimes X_2,\, {\mathbf 1})\\
&\cong& \Hom(X_1\otimes X_2 \otimes V^{*},\, {\mathbf 1})\\
&\cong& \Hom(X_1\boxtimes X_2,\, (V^{**}\boxtimes {\mathbf 1})\otimes A ),
\end{eqnarray*}
therefore, the hexagon identity and Proposition~\ref{equivalence}(a) imply
that 
$$
\iota_V = (u_V\boxtimes \id_{\mathbf 1})\otimes \id_A.
$$

Next, using the definition of $\iota_V$ in \eqref{IV} and 
the naturality of braiding, we compute:
\begin{eqnarray*}
\lefteqn{\iota_{{}^{**}V}^{-1} \circ \iota_{{}^*V}^*  =}\\
&=& ( \rho_{{}^{**}V}^{-1} \circ \Sigma_{{\mathbf 1}\boxtimes {}^{**}V, A}^{-1}  \circ
    (\rho_{{}^{*}V}^*)^{-1} \circ  \Sigma_{A, V \boxtimes {\mathbf 1}}^{-1} ) 
     \circ ( \rho_{{}^{*}V}^* \circ \Sigma_{{\mathbf 1}\boxtimes {}^{**}V, A} \circ
    \rho_V^{**} \circ \Sigma_{A, V ^{**}\boxtimes {\mathbf 1}} )\\
&=& \rho_{{}^{**}V}^{-1} \circ \Sigma_{{\mathbf 1}\boxtimes {}^{**}V, A}^{-1}  \circ
    (\rho_{{}^{*}V}^*)^{-1} \circ \left( \rho_{{}^{*}V}^* \circ 
    \Sigma_{{\mathbf 1}\boxtimes {}^{**}V, A} \circ
    \rho_V^{**} \circ \Sigma_{A, V ^{**}\boxtimes {\mathbf 1}} \right) \circ 
    \Sigma_{A, V ^{**}\boxtimes {\mathbf 1}}^{-1} \\
&=& \rho_{{}^{**}V}^{-1} \circ \rho_V^{**}.
\end{eqnarray*}
On the other hand,
\begin{equation*}
\iota_{{}^{**}V}^{-1} \circ \iota_{{}^*V}^* = 
(u_{{}^{**}V}^{-1} \circ  u_{{}^*V}^* \boxtimes \id_{\mathbf 1})\otimes \id_A,
\end{equation*}
therefore Proposition~\ref{equivalence}(a) and equation \eqref{unimodular delight}
imply the result.
\end{proof}
\end{section}

\begin{section}{Comparison with \cite{EO}}
\label{comp EO}

In this Section we show that our Definition~\ref{dist inv o} 
of a distinguished invertible object $D$ of $\C$ agrees with 
\cite[Definition 2.12]{EO}.

Recall that it was proved in \cite{EO} that in a  finite tensor category 
any projective object is injective and vice versa. In particular the
projective cover $P_0$ of the unit object ${\mathbf 1}\in \C$ coincides with the
injective hull of some object $\tilde D\in \C$. It was shown in \cite{EO}
that $\tilde D$ is an invertible object. 

\begin{theorem} 
\label{comparison}
The object $\tilde D$ is isomorphic to $D$.
\end{theorem}
\begin{proof} Let $I$ be a set indexing the isomorphism classes of simple
objects in $\C$; for $\alpha \in I$ let $L_\alpha, P_\alpha, I_\alpha$ denote
a simple object corresponding to $\alpha$, its projective cover, and its
injective hull. We will assume that $0\in I$ and $L_0={\mathbf 1}$. Let $i$ runs through
$I$. We are going to compute $\dH (P_0\boxtimes L_i, A^*)$ in two ways.
\newline {\em First calculation:}
\begin{eqnarray*}
\dH (P_0\boxtimes L_i,A^*)
&=& \dH (P_0\boxtimes L_i,(D\boxtimes {\mathbf 1})\otimes A)\\
&=& \dH (P_0\otimes L_i,D)
= \left\{ \begin{array}{cc}1&\mbox{if}\; L_i=D,\\
0&\mbox{otherwise}.\end{array}\right. 
\end{eqnarray*}
\newline
{\em Second calculation:}
\begin{eqnarray*}
\dH (P_0\boxtimes L_i,A^*)
&=& \dH (A,{}^*(P_0\boxtimes L_i))\\
&=& \dH ((P_0\boxtimes {\mathbf 1})\otimes A,{\mathbf 1}\boxtimes L_i^*).
\end{eqnarray*}
Let us look closely at the object $(P_0\boxtimes {\mathbf 1})\otimes A$ in $\C
\boxtimes \C^{\op}$.

\begin{lemma} 
\label{preprojective}
The object $(P_0\boxtimes {\mathbf 1})\otimes A$ is injective.
\end{lemma}
\begin{proof} Observe that the functor
$$\Hom (X\boxtimes Y,(P_0\boxtimes {\mathbf 1})\otimes A)=\Hom ((P_0^*\otimes X)\boxtimes
Y,A)=\Hom (P_0^*\otimes X\otimes Y,{\mathbf 1})$$ 
is exact in both variables $X,Y$ since $P_0^*\otimes X\otimes Y$ is injective,
see \cite{EO}. Thus the functor $\Hom(?,(P_0\boxtimes {\mathbf 1})\otimes A)$ is exact,
see \cite{D}. The Lemma is proved.
\end{proof}
\noindent
We continue the proof of the Theorem. By Lemma~\ref{preprojective},
\begin{equation*}
(P_0\boxtimes {\mathbf 1})\otimes A=\sum_{\alpha,\beta \in I}M_{\alpha \beta}
I_\alpha \boxtimes I_\beta
\end{equation*} 
for some non-negative integer
multiplicities $M_{\alpha \beta}$. We have
\begin{eqnarray*}
M_{\alpha \beta}
&=& \dH (L_\alpha \boxtimes L_\beta,(P_0\boxtimes {\mathbf 1})\otimes A)\\
&=& \dH (P_0^*\otimes L_\alpha \otimes L_\beta,{\mathbf 1})\\
&=& \dH (L_\alpha \otimes L_\beta,\,P_0)\\
&=& [L_\alpha \otimes L_\beta :\tilde D],
\end{eqnarray*}
where $[X:L_i]$ denotes the multiplicity of a simple object $L_i$ in the
Jordan-H\"older series of $X$. To calculate $\dH ((P_0\boxtimes {\mathbf 1})\otimes A,
{\mathbf 1}\boxtimes L_i^*)$ it is enough to consider the summands with $I_\alpha =P_0$.
In this case $L_\alpha =\tilde D$ and $[L_\alpha \otimes L_\beta :\tilde D]=
[L_\beta :{\mathbf 1}]$. Thus
\begin{equation}
\dH ((P_0\boxtimes {\mathbf 1})\otimes A,{\mathbf 1}\boxtimes L_i^*)=\left\{ \begin{array}{cc}
1&\mbox{if}\; I_0\; \mbox{covers}\; L_i^*,\\
0&\mbox{otherwise}.\end{array}\right.
\end{equation}
Since $I_0=P_0^*$ covers $\tilde D^{-1}$ the Theorem follows.
\end{proof}

\begin{remark}
Note that Theorems \ref{main} and \ref{comparison} together are equivalent
to \cite[Conjecture~2.15]{EO}. 
\end{remark}

\begin{corollary}
A semisimple  finite tensor category is unimodular.
\end{corollary}

Here is another application of Lemma \ref{preprojective} (which in the case of Hopf
algebras is a well-known statement).

\begin{proposition} 
Let $f:A \to A^{**}$ be a morphism in
$\C\boxtimes\C^{\op}$. Assume that $\Tr (f)\neq 0$. 
Then the category $\cC$ is semisimple.
\end{proposition}
\begin{proof} 
By definition $\Tr (f)$ is the following morphism :
\begin{equation}
\label{trace}
\Tr(f) :  {\mathbf 1}\boxtimes {\mathbf 1} \stackrel{\coev_A}{\longrightarrow} A\otimes A^*
\stackrel{f\otimes \id_A}{\longrightarrow}
A^{**}\otimes A^* \stackrel{\ev_{A^*}}{\longrightarrow} {\mathbf 1}\boxtimes {\mathbf 1}.
\end{equation}
In particular, if $\Tr (f)\ne 0$ then ${\mathbf 1}$ is a direct
summand of $A\otimes A^*$. Hence $P_0\boxtimes {\mathbf 1}$ is a direct summand of
$(P_0\boxtimes {\mathbf 1})\otimes A\otimes A^*$. By Lemma \ref{preprojective}
$(P_0\boxtimes {\mathbf 1})\otimes A$ is projective and therefore 
$(P_0\boxtimes {\mathbf 1})\otimes A\otimes A^*$ is projective. Thus $P_0\boxtimes {\mathbf 1}$
is projective and consequently ${\mathbf 1}$ is projective. Hence $\cC$ is semisimple.
\end{proof}

\end{section}


\begin{section}{Fusion categories}
\label{fusion categories}
In this section we specialize the previous considerations to the case
when the category $\cC$ is semisimple (and thus $\cC$ is a fusion category,
see \cite{ENO}). Let $\{ L_i\}_{i\in I}$ be a set of representatives
of isomorphism classes of simple objects in $\cC$.

Let us describe the structure of the algebra $A$. 
It follows immediately from definitions
that we have a canonical isomorphism 
\begin{equation}
A=\bigoplus_{i\in I}\, L_i\boxtimes {}^*L_i. 
\end{equation}
Using the definitions one describes the multiplication in 
the algebra $A$ in the following way (cf.\ \cite{Mu}): we have 
$A\otimes A=\bigoplus_{i,j\in I}(L_i\otimes L_j)\boxtimes 
{}^*(L_i\otimes L_j)$; for any $m\in I$ the vector
spaces $\Hom (L_i\otimes L_j, L_m)$ and $\Hom ({}^*(L_i\otimes L_j), {}^*L_m)=
\Hom (L_m, L_i\otimes L_j)$ are canonically dual to each other via the
pairing 
$$
\Hom (L_i\otimes L_j, L_m) \otimes \Hom (L_m, L_i\otimes L_j) \to
\Hom (L_m,L_m)=k
$$ 
and hence there is a canonical morphism
$(L_i\otimes L_j)\boxtimes {}^*(L_i\otimes L_j)\to L_m\boxtimes {}^*L_m$.
Then the multiplication $A\otimes A\to A$ is just a direct sum (over $m\in I$)
of all such morphisms.

Next, we describe the canonical isomorphism $A\cong A^{**}$. We have canonically
\begin{equation*}
A^{**} \cong \bigoplus_{i\in I}L_i^{**}\boxtimes {}^{***}L_i \cong
\bigoplus L_i\boxtimes {}^{*****}L_i. 
\end{equation*}
Thus to specify a morphism
$A\to A^{**}$ in $\C\boxtimes \C^{\op}$ is the same as to specify a collection of morphisms
$\psi_i: {}^*L_i\to {}^{*****}L_i$ in $\C$. Recall \cite{Mu, ENO} that 
for any simple object $L$ in $\C$ one defines its {\em squared norm} $|L|^2$ 
as follows: choose an isomorphism $\phi : L\to L^{**}$ (such an isomorphism 
always exists and is unique up to a scaling) and set $|L|^2=\Tr (\phi)
\Tr ((\phi^{-1})^*)$, where $\Tr$ is defined as in equation \eqref{trace}.
Note that $|L|^2$ does not change after rescaling
the evaluation and coevaluation morphisms in $\C$.

\begin{lemma} 
\label{twice}
The canonical isomorphism $A\to A^{**}$ corresponds to the
collection of morphisms $\psi_i: {}^*L_i\to {}^{*****}L_i$ 
characterized by the following property:
for any isomorphism $\phi_i: {}^*L_i\to {}^{***}L_i$ one has 
$\Tr (\phi_i^{-1})\Tr (\phi_i \circ \psi_i^{-1})=|L_i|^2$.
\end{lemma}
\begin{proof} The statement is immediate from definitions since the
isomorphism $A\cong A^{**}$ is the composition of the isomorphism $A\cong A^*$
and of the inverse of the dual of this isomorphism.
\end{proof}

\begin{corollary} 
The trace of the canonical isomorphism $A\cong A^{**}$
is  equal to $\dim (\cC):=\sum_{i\in I}|L_i|^2$.
\end{corollary}

Let us relate the canonical isomorphism $\delta :?^{**}\cong {}^{**}?$
from Theorem~\ref{main}
with squared norms of simple objects of $\C$.

\begin{theorem} 
\label{vitia monstr}
Let $L\in \cC$ be a simple object. The canonical isomorphism
$\delta_L~:~L^{**}\cong {}^{**}L$ can be characterized in the following way:
for any isomorphism $\phi : L^{**}\to L$ one has
$\Tr (\phi^{-1})\Tr (\phi \circ \delta_L^{-1})=|L_i|^2$.
\end{theorem}
\begin{proof} 
Recall that $A$ represents the functor $X\boxtimes Y\mapsto 
\Hom (X\otimes Y,{\mathbf 1})$ and $A^{**}$ represents the functor
$X\boxtimes Y\mapsto \Hom ({}^{**}X\otimes Y^{**},{\mathbf 1})=
\Hom (X\otimes Y^{****},{\mathbf 1})$. It follows immediately from
definitions that the canonical isomorphism $A\to A^{**}$ corresponds
to the natural transformation $X\boxtimes Y\stackrel{\id \boxtimes 
\delta_{Y^{**}}^{-1}}{\longrightarrow}X\boxtimes Y^{****}$. Now the Theorem is
an immediate consequence of Lemma \ref{twice}. 
\end{proof}

\begin{corollary} 
\label{ler}
Let $V\in \cC$ be an object and let $\phi: V\to {}^{**}V$
be a morphism. Then $\Tr (\phi^*)=\Tr (\phi\circ \delta_V^{-1})$. 
\end{corollary}
\begin{proof} It is enough to prove the statement for $V=L_i$ and any
isomorphism $\phi : L_i \to {}^{**}L_i$. But this is an immediate consequence
of Theorem~\ref{vitia monstr}.
\end{proof} 

\begin{remark} 
The following example shows that the statement of Corollary~\ref{ler} is not true if $\C$
is only assumed to be unimodular. 
Let $q$ be a primitive $p$th root of unity and let $U_q(sl_2)$
be the corresponding finite dimensional quantum $sl_2$ Hopf algebra.
Let $H =\mbox{gr}(U_q(sl_2))$ be the associated graded 
Hopf algebra of $U_q(sl_2)$. This is a Hopf algebra of dimension $p^3$ 
defined like $U_q(sl_2)$ except that $EF-FE=0$. Since $U_q(sl_2)$ 
is unimodular, so is $\mbox{gr}(U_q(sl_2))$. 
Then the statement of Corollary~\ref{ler} would say that $\Tr(K)=\Tr(K^{-1})$
in any finite dimensional representation, which is false in $1$-dimensional
representations.
\end{remark}

Recall (see \cite{ENO}) that for any fusion category $\cC$ one constructs
a twice bigger category $\tilde \cC$ called its {\em pivotalization}. By
definition, the objects of $\tilde \cC$ are pairs $(X,f)$ where $X$ is an
object of $\cC$ and $f: X\to X^{**}$ is an isomorphism satisfying 
$f^{**}f=d_{X^{**}}$. It is easy to see that the category $\tilde \cC$ has
a canonical pivotal structure. 

\begin{corollary} The category $\tilde \cC$ is spherical, that is 
$\dim (X)=\dim (X^*)$ for any $X\in \tilde \cC$.
\end{corollary}
\begin{proof} Follows from  Corollary~\ref{ler}.
\end{proof}

\begin{remark} 
The statement of corollary was proved in \cite{ENO} under
assumption that $\dim (\cC)\ne 0$, which is automatically
satisfied in characteristic zero. Thus our result is new only in positive
characteristic.
\end{remark}

In the case when $\mathcal{C}$ is the representation 
category of a semisimple Hopf algebra $H$, 
Corollary~\ref{ler} becomes the following statement, which appears to be new
in the case of positive characteristic (however, see \cite[Theorem 2.2]{LR2}
and \cite[Propositions 9 and 10]{R4} for related statements).

\begin{proposition} 
\label{trtr}
If $H$ is a semisimple Hopf algebra over an algebraically closed field $k$
and $g$ the distinguished grouplike element
in $H$ (so that $g^{-1}xg=S^4(x)$). Then for any element $a\in H$ such that
$axa^{-1}=S^2(x)$ and any irreducible $H$-module $V$, one has 
\begin{equation}
Tr_V(a)=Tr_V(ga).
\end{equation}
\end{proposition}

\begin{remark}
Let us give another, purely Hopf-algebraic proof of this statement. 
Consider a new Hopf algebra $K$ obtained by extending $H$ by a grouplike element 
$a$ such that $a^2=g^{-1}$ and $axa^{-1}=S^2(x)$ \cite{So}. It is well known that 
$S^2$ is an inner automorphism of $H$, so every irreducible 
representation $V$ of $H$ extends to a representation of $K$ (in two ways). 
We must show that $\Tr_V(a)=\Tr_V(ga)$ for either of these extensions. 

To do this, we recall that the left integral 
in $K^*$ is given by the formula \cite[Proposition 2.4(a)]{LR1}
$$
\lambda(x)=\sum_{W\in {\rm Irrep}K}\Tr_W(xa)\Tr_W(a^{-1}).
$$
Also, recall that under right multiplication in $K^*$, 
$\lambda$ changes according to the character $g$. 
So for any $f\in K^*$ we get 
$$
\sum_W \Tr_W(x\1a)f(x\2)\Tr_W(a^{-1})=\sum_W \Tr_W(xa)\Tr_W(a^{-1})f(g). 
$$
(we use the Sweedler notation $\Delta(z)=z\1\otimes z\2$,
implying the summation as usual).
Set $f(z)=\Tr_{V}(za)$, and $x=I$ (the integral of $K$
acting by $1$ in the trivial module; it exists since $K$ is semisimple). 
Then we get 
$$
\sum_W \Tr_W(I\1a)\Tr_V(I\2a)\Tr_W(a^{-1})=\sum_W \Tr_W(Ia)\Tr_W(a^{-1})
\Tr_{V}(ga). 
$$
The right hand side of this equation is obviously equal to $\Tr_V(ga)$
(only the term with $W=k$ survives). As to the left hand side, it can 
be written as 
$$
\sum_W \Tr_{W\otimes V}(Ia)\Tr_W(a^{-1})=
\sum_W \Tr_{W\otimes V}(I)\Tr_W(a^{-1}).
$$
Here the only nonzero summand comes from $W=V^*$, and it yields 
$\Tr_{V^*}(a^{-1})=\Tr_{V^*}(S(a))=\Tr_V(a)$. Thus, $\Tr_V(ga)=\Tr_V(a)$, 
and we are done. 
\end{remark}

\begin{corollary} 
\label{unim} 
If $H$ is a semisimple Hopf algebra over $k$ and 
$S^4=1$ in $H$ then $H^*$ is unimodular. 
\end{corollary}
\begin{proof} In this case $g$ is central, 
so we get $\Tr_V(a)=g_V\Tr_V(a)$, where $g_V$ is the eigenvalue of $g$ on $V$. 
But $\Tr_V(a)$ is nonzero, so $g_V=1$ and hence $g=1$. 
Thus $H^*$ is unimodular. 
\end{proof}

\begin{remark} 
In the case $S^2=1$, this result is contained in \cite{La}.
\end{remark}

\end{section}


\begin{section}{Appendix}

We will use the notation of \cite{EO}. 
Let $\mC$ be a  finite tensor category, and 
$\mM$ be an  abelian category which carries 
the structure of a module category over $\mC$. 
Let $B$ be a finite dimensional algebra over $k$, 
and suppose that an equivalence of 
categories $\mM\to B-{\rm mod}$ is fixed. In this case, 
any object $X\in \mC$ defines an exact functor
$X\otimes: B-{\rm mod}\to B-{\rm mod}$. 
This means that we have a tensor functor 
$F: \mC\to B-{\rm bimod}$, such that $X\otimes M=F(X)\otimes_B M$
for $M\in B-{\rm Mod}$. 

Let $\bar F:=F\circ {\rm Forget}: \mC\to \mbox{Vect}_k$ 
be the fiber functor, and $H:={\rm End} \bar F$. 
Thus $\mC=H-{\rm mod}$ as an abelian category. 
(In fact, $H$ has the structure of a Hopf algebroid,
reflecting the fact that $\mC$ is a  finite tensor category). 
In particular, we have a homomorphism $\mu: B\otimes B^\circ\to
H$, coming from the functor $F$. Thus $H$ is a module over 
$B\otimes B^{\circ}$, where $B^\circ$ is the algebra opposite to $B$,
via $a\circ h:=\mu(a)h$ for all $h\in H$ and 
$a\in B\otimes B^{\circ}$. 

Recall \cite{EO} that $\mM$ is called an {\em exact} module category if
for any projective object $P\in \mC$ and any $M\in \mM$, the
object $P\otimes M$ is projective. Equivalently, $\M$ is exact if
any $\C$-module additive functor from $\M$ is exact.

The main result of this appendix is the following algebraic
characterization of the property of exactness. 

\begin{theorem}
\label{M exact iff}
$\mM$ is exact if and only if $H$ is a projective $B\otimes
B^\circ$-module. 
\end{theorem}
\begin{proof} 
Recall that a module $V$ over a finite dimensional algebra $A$ is
projective if and only if it is flat. Indeed, 
$V\otimes_A W={\rm Hom}_A(V,W^*)^*$ for any two finite
dimensional left $A$-modules $V,W$, and hence 
the functor $V\otimes_A ?$ 
is exact iff so is ${\rm Hom}_A(V,?)$. 
 
Suppose $\mM$ is an exact module category over $\mC$.
We see that to prove the required statement, it suffices to
show that the module $H$ over $B\otimes B^\circ$ is flat.  

Since $H$ is a left $H$-module, we can regard $H$ as an object of
$\mC$. Thus, the functor $H\otimes_{B}\,?$ is exact, and 
for any finite dimensional left $B$-module $M$, the
module $H\otimes_{B}M$ is projective. 

To show that $H$ is flat, let $M$ be a left $B$-module, 
$N$ a right $B$-module, and let us compute the derived tensor
product $H\otimes^L_{B\otimes B^\circ}(M\otimes N)$.
By the K\"unneth formula we have 
$$
H\otimes^L_{B\otimes B^\circ}(M\otimes N)\cong
(H\otimes_B^L M)\otimes^L_{B^\circ}N.
$$  
Now, since the functor $H\otimes_B$ is exact, we have 
$$
(H\otimes_B^L M)\otimes^L_{B^\circ}N\cong
(H\otimes_B M)\otimes^L_{B^\circ}N.
$$
Further, since the module $H\otimes_B M$ is projective, 
it is also flat. Thus,
$$
(H\otimes_B M)\otimes^L_{B^\circ}N
\cong (H\otimes_B M)\otimes_{B^\circ}N.
$$
We conclude that 
$$
H\otimes^L_{B\otimes B^\circ}(M\otimes N)\cong
H\otimes_{B\otimes B^\circ}(M\otimes N),
$$
which implies that $H$ is flat, as desired. 

Conversely, assume that $H$ is a projective module. 
Then $H\cong \oplus_{i,j} V_{ij}\otimes P_i\otimes P_j^\circ$, where 
$V_{ij}$ are vector spaces, and 
$P_i,P_j^\circ$ are the projective covers of irreducible modules 
over $B,\,B^\circ$, respectively. Thus, for any left $B$-module
$M$ we have 
$$
H\otimes_B M=\oplus_{i,} [\oplus_jV_{ij}\otimes
(P_j^\circ\otimes_B M)]\otimes P_i. 
$$
This $B$-module is obviously
projective, so $\mM$ is exact and we are done. 
\end{proof}

\end{section}

\bibliographystyle{ams-alpha}

\end{document}